\documentclass[reqno,11pt]{amsart}
\usepackage{amssymb}

\usepackage[english]{babel}
\usepackage[latin1]{inputenc}
\usepackage[T1]{fontenc}
\usepackage{latexsym,amsfonts}
\usepackage{graphics}
\usepackage{hhline}
\usepackage{mathrsfs}
\usepackage{amsmath}
\usepackage{amsthm}
\usepackage{pifont}
\usepackage{varioref}
\usepackage{textcomp}
\usepackage{lmodern}
\usepackage{euscript}
\usepackage[bookmarks=false]{hyperref}

\makeatletter\@addtoreset{equation}{section}

\newtheorem {theorem}{Theorem}[section]
\newtheorem {definition}[theorem]{Definition}
\newtheorem {remark}[theorem]{Remark}
\newtheorem {proposition}[theorem]{Proposition}
\newtheorem {example}[theorem]{Example}
\newtheorem {lemma}[theorem]{Lemma}
\newtheorem {corollary}[theorem]{Corollary}
\newcommand{\fin}{\hfill $\square$}

\newcommand{\C}{\mathbb C}
\newcommand{\R}{\mathbb R}

\newcommand{\N}{\mathbb N}
\newcommand{\norm}[1]{\left\Vert#1\right\Vert}

\begin{document}
\title[]{A division's theorem on some class of $\mathcal{C}^\infty$-functions}
\author{Mouttaki HLAL}
\address{Département des Sciences économiques,
Faculté des Sciences Juridiques, Economiques et Sociales, BP 145, Mohammadia, Maroc}
\thanks{
AMS classification (2000): 26E10, 32B05.}

\maketitle

 \begin{abstract}
Let $\mathcal{E}_n$ be the ring of the germs of $\mathcal{C}^\infty$-functions
at the origin in $\R^n$. It is well known that if $I$ is an ideal of $\mathcal{E}_n$, generated by a finite number of germs of analytic functions, then $I$ is closed. In this paper we consider an ideal of $\mathcal{E}_n$ generated by a finite number of germs in some class of
$\mathcal{C}^\infty$-functions that are not analytic in $à$, but quasi-analytic and we shall prove that the result holds in this general situation. We remark that the result is not true for a general ideal of finite type of $\mathcal{E}_n$.
\end{abstract}

 \section{Asymptotic expansions and Gevrey asymptotics}

 We denote by $\C[[z]]$ the ring of formal power series with coefficients in $\C$. We say that an analytic function $f$ in a sector $S=\left\{ z\in \C; 0 < |z| < r, \; \theta_0 < \arg z < \theta_1 \right\}$, continuous on $\overline{S}$ admits $\hat f =\sum_{n\in\N} a_n z^n \in \C[[z]]$ as asymptotic expansion at $0$ if for every subsector $S'$ of $S$ there exist $C, M > 0$ such that for every nonnegative integer $n\in \N m $ and every $z\in S'$,
 $$ \left| f(z) - \sum_{p=0}^{n-1} a_p z^p  \right|  \leq CM^n|z|^n.$$
 For $k>0$, $R> 0$, $\eta > 0$, we define the sector
 $$S^{k}_{R,\eta}=\left\{ z\in \C; 0 < |z| < R, \; \mbox{and } \; |\arg z| < \frac \pi{2k} + \eta \right\}.$$
 Let $A_{k,R,\eta}$ be the set of all functions $f(z)$ holomorphic in the sector $S^{k}_{R,\eta}$ continuous on $\overline{S^{k}_{R,\eta}}$ and having  an asymptotic expansion $\hat f(z)$ and $A_k$ the inductive limit of  $(A_{k,R,\eta})_{R,\eta>0}$.

 A function $f \in  A_{k,R,\eta}$  is Gevrey of order $k$ if for all subsector $S'$ of $S^{k}_{R,\eta}$, there exist constants $C_{S'} > 0$, $M_{S'} > 0$ such that $\forall n\in \N$:

 $$\sup_{z\in S'} \left| \frac{f^{(n)}(z)}{n!} \right|  \leq C_{S'} M _{S'}^n (n!)^{\frac 1k},$$
 then we denote by  $\mathcal{G}_{k,R,\eta}$ the algebra of those functions, and by $\mathcal{G}_k$ we denote the inductive limit of $(\mathcal{G}_{k,R,\eta})_{R,\eta>0}$.  $\mathcal{G}_k$  is called the ring of $k$-summable functions in direction $\R_+$. We put $\mathcal{G} = \cup_{k>\frac 12}\mathcal{G}_k$.

 \begin{example} First we will give one way to construct some element of
  $\mathcal{G}$. Let $$S=\left\{ t\in \C; 0 < |t| < r, \; \theta_0 < \arg t < \theta_1 \right\}$$
  be a sector; if $D(0,r) \subset \C$ denote the disk, let $f: S \times D(0,r) \longrightarrow \C$ be holomorphic function. We suppose that:
 $$ \exists c > à, \; A> 0, \; k>0 \quad \mbox{ such that } \forall (t,z) \in  S \times D(0,r), \quad
 \; \; |f(t,z)| \leq c e^{-\frac{A}{|t|^k}}.$$
 Let $h(z)=\int_0^1 \frac{f(t,z)}{t+z}dt$. Then $h$ is holomorphic in $\C\setminus [-1,0]$. If $k>\frac 12$, $h$ is holomorphic in the sector $S^k_{r,\eta}$.
 For each $p\in \N$, we put
 $$ h_p(z)=\int_{\frac 1{(p+1)^{\frac 1k}}}^{\frac 1{p^{\frac 1k}}} \frac{f(t,z)}{t+z}dt.$$
 Then $h_p
 \in \mathcal{H}(S^k)$, where $S^p = S^k_{r,\eta} \bigcup \left\{z\in \C; \; |z| < \frac{r^{\frac 1k}}{(p+1)^{\frac 1k}} \right\}.$ For each $p\in \N$, $h_p$ satisfies the following properties
 \begin{align*}
 & \forall z\in S^p   \quad \norm{h_p}_{S^p}:= \sup_{z\in S^p}|h_p(z)|\leq c \rho^p, \mbox{ where } \; \; c>0, \;\; 0<\rho < 1 \\
 & \forall z\in S  \quad  h(z) = \sum_{p=1}^{\infty}h_p(z).
 \end{align*}
 Following \cite{6}, $h\in \mathcal{G}_{k,r,\eta}$ and the map $h \longrightarrow \hat h \in \C[[z]]$ is injective because the angle of $S^k_{r,\eta}$ is strictly greater than $\frac\pi k$. If $f$ is real then $h|_{]0,r]}$ will be in our algebra $\mathcal{G}$.
 \\
 For example if $f(t,z)= \frac 1 t e^{-\frac 1t}$ then $h|_{]0,r]}$ is analytic in  $]0,r]$  and have a $\mathcal{C}^\infty$-extension to $0$ and its Taylor's expansion at $0$ is a convergent series plus the Euler's series $\sum_{n=0}^{\infty} (-1)^n n! x^{n+1}$.
 \\
 This situation is quite general. Let $r>0$, $0<r<R$, $D_n=\left\{ z\in \C; |z|< \frac 12 (n+1)^{-\frac 1k}\right\}$.
 \end{example}

 \begin{proposition}[\cite{6}] \label{Prop12}
 Let $f: S^k_{R,\eta} \longrightarrow \C$. Suppose that for each $n\in \N$ there exists $f_n: D_n \cup S^k_{R,\eta} \longrightarrow \C$ such that:
 \begin{enumerate}
   \item $f_n$ is holomorphic in $S^n$ and $\norm{f_n}_{S^n} \leq c\rho^n, \quad c>0, \;\; 0<\rho < 1$.
   \item $\forall z\in S^k_{R,\eta} $, $f(z)= \sum_{n\in\N} f_n(z)$.
 \end{enumerate}
 Then $f \in \mathcal{G}_{k,r,\eta}$.
 \end{proposition}

 The converse is also true

 \begin{proposition}[\cite{6}] \label{Prop22}
 Let $(f_n)_{n\in \N}$ be a sequence of holomorphic functions on $D_n \cup S^k_{R,\eta}$ such that
 $\norm{f_n}= \sup_{z\in D_n \cup S^k_{R,\eta}}|f_n(z)|  \leq c\rho^n$
with $c>0$, $0<\rho < 1$, then $f(z)= \sum_{n\in\N} f_n(z)$ is Gevrey of order $k$ in $S^k_{R,\eta} $.
 \end{proposition}

 \begin{lemma}
 The homomorphism $ \mathcal{G}_k \longrightarrow \C[[z]]_k$ is injective.
 \end{lemma}

 \noindent {\bf Proof.} Let $f\in  \mathcal{G}_k$ such that $\hat f =0$, we have that $f\in \mathcal{H}(S^k_{R,\eta})$, and the angle of $S^k_{R,\eta}$ is strictly greater than $\frac \pi k$, so by the result of Gevrey asymptotic functions based on Phragmen-Lindel\"of principle we have $f=0$.
 \fin

 \section*{Gevrey-analytic functions}
  Let $f\in  \mathcal{G}_k$, there exist $R>0$, $0<\rho < 1$, $k>\frac 12$, and $\eta$ small such that $f\in \mathcal{G}_{k,R,\eta}$, by Proposition \ref{Prop12} there exists a sequence of functions $(f_n)_{n\in \N}$, holomorphic in $S^n=D_n \cup S^k_{R,\eta}$ such that
  $\norm{f_n}_{S^n} \rho^{-n} < +\infty$. We put $\underline{\alpha}=(R,\rho,k,\eta)$ and we define
  $$
   \norm{f_n}_{\underline{\alpha}} \rho^{-n} = \inf \left(\sum_{n=0}^\infty \norm{f_n}_{S^n} \rho^{-n} \right)
  $$
  over all representations of $f$ as in Proposition \ref{Prop12}.  For $\mu>0$ and $\underline{\alpha}=(R,\rho,k,\eta)$ as above we define $\mathcal{G}_{\underline{\alpha}}{y_1,\cdots,y_n}_\mu $ as the ring of series $f(z,y)=\sum_{\omega\in\N^n} f_\omega(z)y^\omega$ such that $f_\omega
  \in  \mathcal{G}$ and
  $$
  \norm{f}_{\underline{\alpha},\mu} = \sum_{\omega\in \N^n}\norm{f}_{\underline{\alpha}} \mu^{|\omega|} < +\infty
  $$
  $\left(\mathcal{G}_{\underline{\alpha}}{y_1,\cdots,y_n}_\mu , \norm{}_{\underline{\alpha},\mu}\right)$ is a Banach algebra.
  Then we consider the algebra of Gevrey-analytic functions $\mathcal{G}{y_1,\cdots,y_n}$ as the inductive limit of $\left( \mathcal{G}_{\underline{\alpha}}{y_1,\cdots,y_n}_\mu\right)_{\underline{\alpha}, \mu}$ and we define the inductive limit topology on $\mathcal{G}{y_1,\cdots,y_n}$.

 \begin{proposition} \label{Prop14}
 Let $f\in \mathcal{G}\{y\}$, if $\hat f(0)=0$ then for each $\varepsilon >0$, there exist $\nu,\rho >0$ such that $\norm{f}_{\nu,\rho} \leq \varepsilon$.
 \end{proposition}

  \noindent {\bf Proof.} Let $f\in \mathcal{G}\{y\}$ then $f(z,y)=  \sum_{\omega\in \N^n} f_\omega(z) y^{\omega}$, $f_\omega\in \mathcal{G}$ such that there exist $ \underline{\alpha}, \mu > 0$ such that
  $$
   \sum_{\omega\in \N^n}\norm{f}_{\underline{\alpha}} \mu^{|\omega|} < +\infty .
  $$
  Let $\varepsilon > 0$, then there exist $l\in \N$ such that
  $$
   \sum_{\norm{\omega}>l}\norm{f}_{\underline{\alpha}} \mu^{|\omega|} < \frac\varepsilon 3.
  $$
  Furthermore if we take $\mu$ sufficiently small, we have
    $$
   \sum_{0< \norm{\omega} < l}\norm{f}_{\underline{\alpha}} \mu^{|\omega|} < \frac\varepsilon 3.
  $$
  In the other hand, since $\hat f_0(0)=0$, we can suppose that $f_0(z)= \sum_{n=0}^\infty f_{0,n}(z)$ with $f_{0,n}(0)=0$ for all $n\in \N$ (since
  $f_0(z)= \sum_{n=0}^\infty [f_{0,n}(z)-f_{0,n}(0)]$ and $f_0(0)=0$) such that $\sum_{n=0}^\infty \norm{f_{0,n}}_{S^n} \rho^{-n} <+\infty$. Then there exists $N_0\in \N$ such that $\sum_{n> N_0}^\infty \norm{f_{0,n}}_{S^n} \rho^{-n} <\frac \varepsilon 6$.  Since $f_{0,n}(0)=0$, then if we make $R$ small, we have that $ \norm{f_{0,n}}_{S^n} \rho^{-n} <\frac \varepsilon {6(N_0+1)}$ for $n=0, \cdots, N_0$. Hence,  $\norm{f_{0}}_{\underline{\alpha}}  <\frac \varepsilon 3$. Finally, we have that $\norm{f}_{\underline{\alpha},\mu}  < \varepsilon $.
  \fin

  \begin{corollary}
  Let $f\in \mathcal{G}\{y\}$, if $\hat f(0) \ne 0$. Then there exists $g \in \mathcal{G}\{y\}$ such the $fg=1$.
  \end{corollary}

    \noindent {\bf Proof.}  Put $\hat f(0) = a_0 \ne 0$, we define $\varphi(\xi)=\frac 1{\xi+a_0}$, then we have $\varphi\in \mathcal{H}(D(0,r))$, with $r<|a_0|$. We put $g=f-a_0 \in \mathcal{G}\{y\}$. We have $\hat g(0)=0$. By Proposition \ref{Prop14}, there exist $\nu, \rho >0$ such that $\norm{g}{\nu,\rho} \leq \frac{a_0}2$, then $\varphi(g) \in \mathcal{G}\{y\}$, i.e., $\frac 1 f \in \mathcal{G}\{y\}$.
    \fin

  \section{Algebraic properties of $\mathcal{G}\{y\}$}
  \subsection{Weierstrass preparation theorem}

   Let $y=(y_1,\cdots,y_n)$ and $f\in \mathcal{G}\{y\}$. We say that $f$ is regular of order $p$ in $y_n$ if $\hat f(0,0,y_n)$ is regular of order $p$ in $y_n$.

  \begin{theorem}\label{Thm21}
  If $f\in \mathcal{G}\{y\}$ is regular of order $p$ in $y_n$, then for every $\varphi \in \mathcal{G}\{y\}$, there exist $Q\in \mathcal{G}\{y\}$, $R \in \mathcal{G}\{y_1,\cdots,y_{n-1}\}[y_n]$, $\deg_{y_n} R <p$ such that $\varphi = fQ+ R$. Furthermore $Q$ and $R$ are uniquely determined.
  \end{theorem}

  \noindent {\bf Proof.}
  Let $f\in \mathcal{G}\{y\}$ be a regular of order $p$ in $y_n$. We can write $f= \sum_{m=0}^\infty f_m(z,y_1,\cdots,y_{n-1}) y_n^m$, by hypothesis we have $f_m(0)=0$ for $m<p$ and $f_p(0)\ne 0$. After dividing $f$ by a unit element in $\mathcal{G}\{y\}$, we can suppose that $f_p(z,y_1,\cdots,y_{n-1})\equiv 1$  and $f\in \mathcal{G}_{\underline{\alpha}}\{y\}_\mu$,
  $$ \mathcal{G}_{\underline{\alpha}}\{y\}_\mu = \left\{ \varphi(z,y)=  \sum_{\omega\in \N^n} \varphi_\omega(z) y^{\omega}, \quad \varphi_\omega\in \mathcal{G}/ \quad  \norm{\varphi} =    \sum_{\omega\in \N^n}\norm{\varphi_\omega}_{\underline{\alpha}} \mu^{|\omega|} < +\infty \right\} .$$
  Let $\delta >0$, $\delta < \mu$ such that $\norm{f_m}_{\underline{\alpha}}\leq M_1 \delta$, for $m<p$ (since $f_m(0)=0$ for $m <p$), $M_1$ is a constant independent of $\delta$ and $\mu$, and
  $$ \norm{ \sum_{m>p}^\infty f_m(z,y_1,\cdots,y_{n-1}) y_n^m } \leq M_2\mu^{p+1}.$$
  Hence $\norm{f-y_n^p} = \norm{f - f_p(z,y_1,\cdots,y_{n-1}) y_n^p} \leq M_3 \delta + M_2 \mu^{p+1}$ ($M_3$, $M_2$ are constants independent of $\delta$ and $\mu$).
  \\
  Given $\varepsilon >0$, $\varepsilon < 1$, we can choose $\delta =\delta(\mu)$ such that
  $$ M_2 \mu^{p+1} + M_3 \delta < \varepsilon \mu^{p}.$$
  Such $\delta$ exists if $\mu$ is sufficiently small. Thus
 $$  \norm{f-y_n^p} < \varepsilon \mu^{p}.$$
 On the other hand, if $\varphi \in \mathcal{G}_{\underline{\alpha}}\{y\}_\mu$, we write
 $$ \varphi =Q(\varphi) y_n^p + R(\varphi),$$
 where $Q(\varphi), R(\varphi) \in \mathcal{G}_{\underline{\alpha}}\{y\}_\mu$ and $R(\varphi)$ is a polynomial in $y_n$ of degree $<p$, then
  $$ \norm{\varphi} =  \norm{Q(\varphi)} \mu^p +  \norm{R(\varphi)} . $$
  Let us consider the linear operator
 \begin{align*}
  D &:  \mathcal{G}_{\underline{\alpha}}\{y\}_\mu \longrightarrow  \mathcal{G}_{\underline{\alpha}}\{y\}_\mu\\
    & \varphi \longmapsto D\varphi = \norm{Q(\varphi)} f  +  \norm{R(\varphi)}.
 \end{align*}
 Then
 \begin{align*}
  \norm{D\varphi - \varphi}  &= \norm{Q(\varphi) (f - y_n^p)} \\
            & \leq \varepsilon \norm{Q(\varphi)} \mu^p \leq \varepsilon \norm{\varphi} .
 \end{align*}
 So that we have
 $$ \norm{D - id_{\mathcal{G}_{\underline{\alpha}}\{y\}_\mu }} <1 .$$
 Since $\mathcal{G}_{\underline{\alpha}}\{y\}_\mu $ is a Banach space,  $D$ is invertible, in other words if $\varphi \in \mathcal{G}_{\underline{\alpha}}\{y\}_\mu$ there exists $\psi \in \mathcal{G}_{\underline{\alpha}}\{y\}_\mu $ such that
 $$ \varphi= D\psi =Q(\psi) f+ R(\psi).$$
  \fin

    \begin{remark}
 Let $\mathcal{G}^R\{y_1,\cdots,y_n\}= \{ f\in \mathcal{G}\{y_1,\cdots,y_n\}; \mbox{ such that } \hat f \in \R[[z,y_1, \cdots , y_n]] \}$. Then the last theorem holds for $\mathcal{G}^R\{y_1,\cdots,y_n\}$.
    \end{remark}

    As a consequence of Theorem \ref{Thm21}, we have
     \begin{proposition} \label{Prop22}
   $\mathcal{G}\{y\}$ is a noetherian, local and regular ring of dimension $n+1$.
     \end{proposition}
 We have also
     \begin{proposition} \label{Prop23}
$\C[[z,y]]$ is the completion of the ring $\mathcal{G}\{y\}$ for the $\underline{m}$-adic topology, where $\underline{m}=(z,y_1,\cdots,y_n)$.
     \end{proposition}

   \begin{corollary}
   $\mathcal{G}\{y\}$ is a normal ring.
  \end{corollary}

  \noindent {\bf Proof.}
  $\mathcal{G}\{y\}$ is a local and regular ring, then is factorial (see \cite{4}), so that since $\mathcal{G}\{y\}$ is a noetherian and factorial, then is normal.
  \fin \\
  Since $\mathcal{G}\{y\}$ is a noetherian and local ring, we have
  \begin{corollary}
 $\C[[z,y]]$  is faithfully flat over $\mathcal{G}\{y\}$.
  \end{corollary}

\begin{proposition} \label{Prop26}
$\mathcal{G}\{y\}$ is a henselian ring.
\end{proposition}

  \noindent {\bf Proof.}  Let $f\in \mathcal{G}\{y\}$, $f(u)=u^p + \sum_{i=1}^p a_i$ and $x\in \mathcal{G}\{y\}/\underline{m}\sim \C$ such that $\overline{f(x)}=0$ and $\frac{\partial \overline f}{\partial u}(x) \ne 0$. If $\alpha \in \mathcal{G}\{y\}$, with $\overline{\alpha}=x$ and $v=\alpha -u$, and $g(v)=f(u)$ we have $\overline g(0)=0$ and $\frac{\partial \overline g}{\partial v}(x) \ne 0$, then $g\in \mathcal{G}\{y\}$ is regular of order $1$ in $v$, so that there exist $R \in \underline{m}$ and $Q$ invertible in $\mathcal{G}\{y,v\}$ such that $g=Q(v-R)$, so we have$f=\widetilde{Q}(u-v+R)$, then if $a=\alpha +R$, we obtain $f(a)=0$ and $\overline a =x$.
 \fin

  \begin{definition}
We say that $A$ is pseudo-geometric ring if $A$ is noetherian and if, for every prime ideal $\mathcal {P}$ of $A$ , $A/\mathcal{P}$ satisfies the finiteness condition for integral extensions, i.e., if for every integral extension of $A/\mathcal{P}$ such that its quotient field id finite over the quotient field of
  $A/\mathcal{P}$ then is finite over $A/\mathcal{P}$.
  \end{definition}

  \begin{theorem}\label{Thm28}
  $\mathcal{G}\{y\}$ is pseudo-geometric ring.
  \end{theorem}

 \noindent {\bf Proof.} Let $k \in \{0, \cdots , n\}$ and $\mathcal{P}$ be a prime ideal of height $k+1$, then we have
 $ \mathcal{G}\{y\}/ \mathcal{P}$ is a finite over $\mathcal{G}\{y_1, \cdots,y_{n-k}\}$.
 \\
  In what follows we denote by $y'$ the $(n-k)$-uplet $(y_1, \cdots,y_{n-k})$. If $[B]$ (respectively $[\mathcal{G}\{y'\}]$) denote the quotient
  field of $B$ (respectively $\mathcal{G}\{y'\}$), $[B]$ is a finite algebraic extension of $[\mathcal{G}\{y'\}]$. It follows from theorem on primitive element that
  $$
   [B] = [\mathcal{G}\{y'\}][x],
   $$
  where $x=\frac{b'}{b}$, $(b,b') \in B^2$ and $b\ne 0$. If $b^{"}\in B\setminus\{0\}$ such that $bb^{"}\in \mathcal{G}\{y'\}$, then generate $[B]$ over $[\mathcal{G}\{y'\}]$, then we may suppose that $x\in B$. Let  $P$ be the minimal polynomial of $x$, $\Delta$ its determinant, and let $s=\deg P [[B]:[\mathcal{G}\{y'\}]]$.
  $$
  P(x)= x^s + \sum_{i=0}^{s-1} a_i x^i .
  $$
  Let $\sigma_1, \cdots , \sigma_s$ be the $[\mathcal{G}\{y'\}]$-isomorphisms of $[B]$ in an algebraic closure of $[\mathcal{G}\{y'\}]$, since $x\in B$ then $x$ is integral over $\mathcal{G}\{y'\}$, and also the symmetric functions $a_i$ of $\sigma_i(x)$, therefore $a_j\in \mathcal{G}\{y'\}$ (since $a_j\in [\mathcal{G}\{y'\}]$ and $\mathcal{G}\{y'\}$ is normal). We have
  $$
  \Delta = \prod_{i<j} (\sigma_i(x)-\sigma_j(x))^2 = \det[\sigma_i(x)]^2 .
  $$
  Then $\Delta \ne 0$ and belongs to $\mathcal{G}\{y'\}$. Let $\alpha\in B$. There are elements $(b_j)_{0\leq j \leq s-1}$ of $\mathcal{G}\{y'\}$ such that
  $$
  \alpha = \sum_{j=0}^{s-1} b_j x^j .
  $$
  Then, we have $\sigma_i(x) =  \sum_{j=0}^{s-1} \sigma_j(x) x^j$ for $1 \leq i \leq s$. If we solve this system for $b_j$, we conclude that $\Delta b_j$ ($0\leq j \leq s-1$) is integral over $\mathcal{G}\{y'\}$, so that $\Delta b_j\in \mathcal{G}\{y'\}$ (since $\Delta b_j\in [\mathcal{G}\{y'\}]$ and $ \mathcal{G}\{y'\}$ is normal), therefore we have
  $$\Delta B \subset  \mathcal{G}\{y'\} x^{s-1} + \cdots + \mathcal{G}\{y'\} x. $$
  Then $\Delta B$ and so that $B$ is finite over $\mathcal{G}\{y'\}$, i.e., $\mathcal{G}\{y'\}$ is pseudo-geometric.
 \fin

  \begin{corollary}\label{Cor29}
   Let $\mathcal{P}$ be a prime ideal of $\mathcal{G}\{y\}$. If $\widehat{\mathcal{P}}$ denote the ideal of $\C[[z,y]]$ generated by $\mathcal{P}$, then $\widehat{\mathcal{P}}$ is also prime. Furthermore $\mathcal{P}$ and  $\widehat{\mathcal{P}}$ have the same height.
  \end{corollary}

    \noindent {\bf Proof.} Since $\mathcal{G}\{y\}$ is henselian and pseudo-geometric, then $\mathcal{G}\{y\}$ is a Weierstrass ring, so that $ \widehat{\mathcal{P}}$ is also prime (see \cite{3}). Furthermore since $\C[[z,y]]$ is faithfully flat over $\mathcal{G}\{y\}$, then $\widehat{\mathcal{P}}\cap \mathcal{G}\{y\}= \mathcal{P}$, so that $\mathcal{P}$ and  $\widehat{\mathcal{P}}$ have the same height.
    \fin

  \subsection{Artin's approximation theorem}
    As the classical case we have Artin's theorem.

   \begin{theorem}[Artin in Gevrey]\label{Thm22}
Let $x=(x_1, \cdots , y_n)$ and $y=(y_1, \cdots , y_p)$. Let $f_1(z,x,y)$, $\cdots $, $f_q(z,x,y)$ be in $\mathcal{G}\{x,y\}$. Suppose there exist $\widetilde{y}_1(z,x), \cdots , \widetilde{y}_p(z,x) \in \R[[z,x]]$ such that
$$f_1(z,x,\widetilde{y}(z,x))= \cdots = f_q(z,x,\widetilde{y}(z,x)) =0,$$
where $\widetilde{y}(z,x)=(\widetilde{y}_1(z,x), \cdots , \widetilde{y}_p(z,x) )$. Then for each $\nu\in\N$ there exists
$(y_1(z,x), \cdots , y_p(z,x) ) \in \mathcal{G}\{x,y\}^p$ such that
$$f_1(z,x,y(z,x))= \cdots = f_q(z,x,y(z,x)) =0$$ and $y_j - \widetilde{y}_j \in \underline{m}^{\nu+1}$.
  \end{theorem}

\noindent {\bf Proof.}
The proof is similar to the analytic case (see \cite{2}).
    \fin

    \section{Fréchet modules on $\mathcal{E}(\Omega)$.}
    Let us recall the following result which we will use at follows.

      \begin{proposition} \label{Prop31}
   Let $\mathcal{P}$ be a prime ideal in $\mathcal{G}\{y\}$ of height $k$, and $f_1, \cdots , f_k \in \mathcal{P}  $.
   Then the maximal ideal $\mathcal{P}\mathcal{G}\{y\}_{\mathcal{P}}$ in $\mathcal{G}\{y\}$ is generated by $f_1, \cdots , f_k$ if and only if there exists a jacobian $\frac{D(f_1, \cdots , f_k)}{D(w_1, \cdots , w_k)} \notin \mathcal{P}$, where $w_i$ is $z$ or $y_i$ for $ 1\leq i\leq n$. In this case, there exists
   $\delta\in  (\mathcal{G}\{y\}\setminus \mathcal{P}) \cap \mathcal{J}(f_1, \cdots , f_k)$ such that
   $$\delta \mathcal{P} \subset   < f_1, \cdots , f_k> ,$$
   where $\mathcal{J}(f_1, \cdots , f_k)$ is the ideal generated over $ \mathcal{G}\{y\}$ by $(f_1, \cdots , f_k)$ and all of the jacobians
   $\frac{D(f_1, \cdots , f_k)}{D(w_1, \cdots , w_k)}$ for  $ 1\leq i_1 < i_2 \cdots <i_k \leq n+1 $ and $w_i$ is $z$ or $y_i$ for $ 1\leq i\leq n$.
     \end{proposition}

  \noindent {\bf Proof.}
The proof is similar to the analytic case (see \cite{5}).
    \fin \\
  Let $\Omega$ be an open set in $\R^+\times \R^n$, if we denote be $\mathcal{E}$ the sheaf of real $\mathcal{C}^\infty$-functions
   on $\Omega$, and by $\mathcal{H}$ the sheaf of real analytic functions on $\Omega$. We define the subsheaf $\underline{\mathcal{G}}$ of $\mathcal{E}$
   by: For $(x_0,y_0) $ such that $x_0 > 0$, we put $\underline{\mathcal{G}}_{x_0,y_0}=\mathcal{H}_{x_0,y_0}$ and $\underline{\mathcal{G}}_{0,y_0}=     \mathcal{G}\{y-y_0\}$.
   \\
   Using the same notations as Malgrange \cite{2}, we define the sheaf $\widetilde{\mathcal{G}}$ as follows: For any subset of $\Omega$, with $U=[0,\varepsilon[ \times U'$, where $U'$ is an open neighborhood of $0$ in $\R^n$, we associate
   $$
   \widetilde{\mathcal{G}}(U) = \prod_{(z,y)\in U} \underline{\mathcal{G}}_{(z,y)} .
    $$
  Let $(z,y)\in \Omega$, we denote by $\mathcal{F}_{z,y}$ the completion of $\underline{\mathcal{G}}_{(z,y)}$ for the Krull topology:
  $\mathcal{F}_{z,y}=\R[[z,y]]$, and we define the sheaf $ \widetilde{\mathcal{F}} = \prod_{(z,y)\in U}\mathcal{F}_{z,y}$.
  \\
  We have a coherence theorem analogous to the Oka's coherence theorem for the analytic case.

  \begin{theorem}\label{Thm32}
The sheaf $\widetilde{\mathcal{G}}$ is flat over $\underline{\mathcal{G}}$, i.e., for any $(z,y)\in \Omega$, the module $\widetilde{\mathcal{G}}_{(z,y)}$ is flat over $\underline{\mathcal{G}}_{(z,y)}$.
  \end{theorem}

In the Proof we will consider only the point $(0,y_0)$, for the other points we have the analytic case.\\

    \noindent {\bf Proof.}
    Let $U$ be an open neighborhood of $(0,y_0)$. Let $\underline{f} = (f_1, \cdots , f_k)$ and $\underline{f}_{(0,y_0)} = (f_{1,(0,y_0)}, \cdots , f_{k,(0,y_0)}) \in \mathcal{G}\{y-y_0\}$ such  that $f_{1,(0,y_0)}, \cdots , f_{k,(0,y_0)}$ be a $\mathcal{G}\{y-y_0\}$-sequence of $\mathcal{G}\{y-y_0\}$.
    We shall prove that
    $$ Tor_1^{\mathcal{G}\{y-y_0\}}  \left( \mathcal{G}\{y-y_0\} / \underline{f}_{(0,y_0)} , \widetilde{\mathcal{G}}_{(0,y_0)} \right) =0.$$
   Let $\mathcal{I}$ be the ideal generated in $\mathcal{G}\{y-y_0\} $ by $f_1, \cdots , f_k$. By using Theorem II. 5.3. in \cite{5} for the formal case and Corollary \ref{Cor29}, and reducing $U$ if necessary, we have for any $(z,y)\in U$,
   $$ht(\mathcal{I}_{(z,y)}) \geq ht(\mathcal{I}_{(0,y_0)}) =k .$$
 Then we have either $\mathcal{I}_{(z,y)} = \mathcal{G}\{y-y_0\} $  or $f_{1,(0,y_0)}, \cdots , f_{k,(0,y_0)}$ be a system of parameters of $\mathcal{G}\{y-y_0\} $. Consequently, in all cases we have that the module $\mathcal{R}_{\mathcal{G}\{y-y_0\}} (\underline{f}_{(0,y_0)})$ relations of $\underline{f}_{(0,y_0)}$ in
 $\mathcal{G}\{y-y_0\} $ is generated by the obvious relations, then we have that
 $$ \mathcal{R}_{\mathcal{G}(0,y_0)} \left( \underline{f}_{(0,y_0)} \right) = \mathcal{R}_{\mathcal{G}\{y-y_0\}}   \left( \underline{f}_{(0,y_0)} \right) \widetilde{\mathcal{G}}_{(0,y_0)}.$$
 Hence
 $$
 Tor_1^{\mathcal{G}\{y-y_0\}}  \left( \mathcal{G}\{y-y_0\} / \underline{f}_{(0,y_0)} , \widetilde{\mathcal{G}}_{(0,y_0)} \right)  =
 \mathcal{R}_{\mathcal{G}(0,y_0)} \left( \underline{f}_{(0,y_0)} \right) / \mathcal{R}_{\mathcal{G}\{y-y_0\}}   \left( \underline{f}_{(0,y_0)} \right) \widetilde{\mathcal{G}}_{(0,y_0)} =0.
 $$

    \begin{remark}\label{Rem1}
Since every $\underline{\mathcal{G}}_{(z,y)}$-sequence of $\underline{\mathcal{G}}_{(z,y)}$  is $\mathcal{G}_{(z,y)}$-sequence of $\mathcal{F}_{(z,y)}$, we prove a similar result as Theorem \ref{Thm32} for $ \widetilde{\mathcal{F}}$, i.e., $ \widetilde{\mathcal{F}}$  is flat over $\underline{\mathcal{G}}$.
    \end{remark}

    \subsection{Fréchet's module}

      \begin{definition}
Let $M$ be an $\mathcal{E}(\Omega)$-module. $M$ is a Fréchet module if $M$ is of finite presentation on  $\mathcal{E}(\Omega)$ and the map
$$ i: M \longrightarrow M \bigotimes_{\mathcal{E}(\Omega)}   \widetilde{\mathcal{F}}(\Omega)$$ is injective.
     \end{definition}

    \subsection{Local Fréchet's modules }

      \begin{definition}
Let $M$ be an $\mathcal{E}_a$-module ($a\in \Omega$). $M$ is a Fréchet module if $M$ is of finite presentation on  $\mathcal{E}_a$ and the map
$$ i: M \longrightarrow M \bigotimes_{\mathcal{E}_a}   \widetilde{\mathcal{F}}_a$$ is injective.
     \end{definition}

      \begin{proposition}[\cite{5}] \label{Prop322}
      Let $ 0  \longrightarrow M^{'} \longrightarrow M \longrightarrow M^{"} \longrightarrow 0$ be a exact sequence of $\mathcal{E}_a$-modules. If
      $ M^{'}$ and $M^{"}$ are Fréchet modules, then $M$ is a Fréchet module.
     \end{proposition}

     \subsection{Fréchet's modules of finite real dimension}
Let $\underline{m}_a$ be the maximal ideal of $\mathcal{E}_a $. $\underline{m}_a$ is the ideal of germs of all functions of $\mathcal{E}_a $ which
vanish at $a$, then $\underline{m}_a$ is generated by $y_1-a_1, \cdots ,y_{n+1}-a_{n+1}$ ($a=(a_1, \cdots ,a_{n+1})$), $\mathcal{E}_a/ \underline{m}_a$ is a Fréchet module.
\\
Let $N$ be an $\mathcal{E}_a$-module, and suppose that $\dim_R(N)<+\infty$. The decreasing sequence $\underline{m}_a^iN$ is stationary, then there exists $i$ such that $\underline{m}_a^iN= \underline{m}_a \cdot \underline{m}_a^i\cdot N$. By Nakayama's lemma, we have $\underline{m}_a N=0$. Hence $N$ is of finite length over the noetherian ring $\mathcal{F}_a$. Then there exists an increasing sequence $N_0, \cdots , N_s$ of submodules of $N$ such that $N_0=0$, $N_s=N$ and for $i\in [0,s-1]$, $N_{i+1}/N_{i} \sim \mathcal{E}_a/ \underline{m}_a$, then $N$ is a Fréchet module (Proposition \ref{Prop322}).
\\
Let $\mathcal{E}_{n+1}$ denote the ring of germs at $0\in\R^{n+1}$ of $\mathcal{C}^\infty$-functions with real values.

  \begin{theorem}\label{Thm331}
Let $M$ be a finite $\mathcal{G}\{y\}$-module, then $M \bigotimes_{\mathcal{G}\{y\}} \mathcal{E}_{n+1}$ is a Fréchet module over $\mathcal{E}_{n+1}$.
  \end{theorem}

    \noindent {\bf Proof.}
    Let  $\widetilde{\mathcal{F}}_{n+1}$ be the ring of germs at $0\in\R^{n+1}$ of collections of formal power series at each point near $
    0$. The module $M \bigotimes_{\mathcal{G}\{y\}} \mathcal{E}_{n+1}$ is of finite presentation over $\mathcal{E}_{n+1}$ (since $\mathcal{G}\{y\}$ is noetherian, then $M$ is of finite representation over $\mathcal{G}\{y\}$). We will prove that the map
    $$
    i: M \bigotimes_{\mathcal{G}\{y\}} \mathcal{E}_{n+1} \longrightarrow M \bigotimes_{\mathcal{G}\{y\}} \widetilde{\mathcal{F}}_{n+1}
    $$
    is injective. The sequence $ 0  \longrightarrow \mathcal{E}_{n+1} \longrightarrow \widetilde{\mathcal{F}}_{n+1} \longrightarrow  \widetilde{\mathcal{F}}_{n+1}/\mathcal{E}_{n+1} \longrightarrow 0$ is exact, then we have an exact sequence:
    $$
    Tor_1^{\mathcal{G}\{y\}} (M, \widetilde{\mathcal{F}}_{n+1}) \longrightarrow Tor_1^{\mathcal{G}\{y\}}(M,\widetilde{\mathcal{F}}_{n+1}/\mathcal{E}_{n+1})
     \longrightarrow  M \bigotimes_{\mathcal{G}\{y\}} \mathcal{E}_{n+1} \longrightarrow M \bigotimes_{\mathcal{G}\{y\}} \widetilde{\mathcal{F}}_{n+1}.
     $$
    Since we have $ Tor_1^{\mathcal{G}\{y\}} (M, \widetilde{\mathcal{F}}_{n+1})=0$ (Remark \ref{Rem1}), then $ M \bigotimes_{\mathcal{G}\{y\}} \mathcal{E}_{n+1}$ is a Fréchet module if and only if $Tor_1^{\mathcal{G}\{y\}}(M,\widetilde{\mathcal{F}}_{n+1}/\mathcal{E}_{n+1})=0$. We will prove this result by induction on $\dim(M)$.  If $\dim(M)=0$, then $ M \bigotimes_{\mathcal{G}\{y\}} \mathcal{E}_{n+1}$ is an $\R$-vector space of finite dimension, then $ M \bigotimes_{\mathcal{G}\{y\}} \mathcal{E}_{n+1}$ is a Fréchet module (see 3.3).
    \\
    Suppose now that $\dim(M) =n-k>0$. Then after using the exactness of sequence of $Tor$ we can suppose that $M=\mathcal{G}\{y\}/\mathcal{P}$, where $\mathcal{P}$ is a prime ideal of height $k+1$.
    \\
    Let $\varphi_1, \cdots ,\varphi_s$ Gevrey-analytic functions on $\Omega$ which generate $\mathcal{P}$ over $\mathcal{G}\{y\}$. Then using Proposition 2 for suitable $\varphi_i$, we suppose that there exists $\delta \in \mathcal{G}\{y\}$ such that
    $$
    \delta\varphi_j \in <\varphi_1, \cdots ,\varphi_{k+1} >  \quad \mbox{for } j=k+2, \cdots, s, \quad \delta \in \mathcal{J}_{k+1}(\varphi_1, \cdots ,\varphi_{k+1}).
    $$
    Hence the germ $\delta_0$ of $\delta$ doesn't belong to $\mathcal{P}$. Since $\delta$ is not a zero divisor in $\mathcal{G}\{y\}/\mathcal{P}$, then
    $\delta$ is not a zero divisor in $\widetilde{\mathcal{F}}_{n+1}/ \mathcal{P}\widetilde{\mathcal{F}}_{n+1}$ ($\widetilde{\mathcal{F}}_{n+1}$ is flat over
    $\mathcal{G}\{y\}$, Remark \ref{Rem1}). In other words, by making $\Omega$ smaller, we have $\forall (z,y)\in V(I)$, $\delta_{z,y} =T_{z,y}\delta$
     is not a zero divisor in $\mathcal{F}_{z,y}/T_{z,y} I$, where $V(I)$ is the set of zeros of the ideal $I$. Putting $I'=I+ \delta \cdot \mathcal{E}(\Omega)$,
     if $\Omega$ is sufficiently small; $I'$ is closed (since the height of $\mathcal{P} + \delta \cdot \mathcal{G}\{y\} > k+1$, and we apply the induction hypothesis), hence $I$ is closed (see V. 5.6. \cite{5}), then the ideal  $\mathcal{P}\cdot \mathcal{E}_{n+1}$ induced by $I$ at the origin is closed, therefore $ M \bigotimes_{\mathcal{G}\{y\}} \mathcal{E}_{n+1}= \mathcal{E}_{n+1}/ \mathcal{P}\cdot \mathcal{E}_{n+1}$ is a Fréchet module.
     \fin

    \begin{corollary}
    Let $\mathcal{I}$ be a submodule of $\mathcal{G}\{y\}^p$, then $\mathcal{I}\cdot \mathcal{E}_{n+1}$ is closed in $\mathcal{E}_{n+1}^p$.
  \end{corollary}

    \noindent {\bf Proof.}
    If we put $M=\mathcal{G}\{y\}^p/\mathcal{I}$, then $M$ is a finite $\mathcal{G}\{y\}^p$-module. Hence by Theorem \ref{Thm331},
    we have $ M \bigotimes_{\mathcal{G}\{y\}} \mathcal{E}_{n+1}$ is a Fréchet module, but  $ M \bigotimes_{\mathcal{G}\{y\}} \mathcal{E}_{n+1} =
    \mathcal{E}_{n+1}^p/ \mathcal{I}\cdot \mathcal{E}_{n+1}$, then $\mathcal{I}\cdot \mathcal{E}_{n+1}$ is closed in $ \mathcal{E}_{n+1}^p$.
    \fin \\

Let $f_1,  \cdots, f_q$ be in $\mathcal{G}\{y\}$, then there exist $\varepsilon >0$, $U$ an open neighborhood of $0$ in $\R^n$ such that $f_1,  \cdots, f_q$ are Gevrey-analytic in  $[0,\varepsilon[\times U$.

  \begin{corollary} Let $\Omega$ be an open set in $\R_+\times \R^n$ such that $\Omega=[0,\varepsilon[\times U$, $\varepsilon >0$, where $U$ is an open neighborhood of $0$, and $f_1,  \cdots, f_q$ are Gevrey-analytic in $\Omega$. Let $\Phi \in \mathcal{E}(\Omega)$, then $\Phi$ can be written in the form $\Phi = \sum_{i=1}^q f_i \Psi_i$, where $\Psi_i \in \mathcal{E}(\Omega)$, if and only if for all $(t,a)\in \Omega$, the Taylor expansion $T_{(t,a)}f$ of $f$ at $(t,a)$ belongs to the ideal generated in $\mathcal{F}_{(t,a)}$ by $T_{(t,a)}f_1, \cdots , T_{(t,a)}f_q$.

  \end{corollary}

 \end{document}